\documentclass{sig-alt-gov2}
\usepackage{amsfonts,latexsym,amssymb,amsmath,epsfig}
\newlength{\jmr}
\newlength{\barvinok}
\newtheorem{dfn}{Definition}[section]
\newtheorem{brent}{Brent-Salamin Theorem}

\newtheorem{rem}[dfn]{Remark}
\newtheorem{prop}[dfn]{Proposition}
\newtheorem{thm}[dfn]{Theorem} 

\newtheorem{cor}[dfn]{Corollary}
\newtheorem{ex}[dfn]{Example}

\newtheorem{algor}[dfn]{Algorithm}

\newlength{\smale}
\newcommand{\thth}{^{\text{\underline{th}}}}

\newcommand{\hpt}{{\mathbf{HPTAS}_\R}}
\newcommand{\hhpt}{{\mathbf{HPTAS}}}
\newcommand{\fpt}{{\mathbf{FPTAS}}}
\newcommand{\np}{{\mathbf{NP}}}

\newcommand{\opt}{{\mathbf{SUP}}}
\newcommand{\feas}{{\mathbf{FEAS}}}
\newcommand{\oopt}{{\mathbf{FSUP}}}

\newcommand{\gln}{\mathbb{G}\mathbb{L}_n}

\newcommand{\bp}{{\mathbf{BP}}}
\newcommand{\ch}{{\mathbf{CH}}}
\newcommand{\ppoly}{\pp\mathbf{/Poly}} 
\newcommand{\pspoly}{\pspa\mathbf{/Poly}} 
\newcommand{\npoly}{\np\mathbf{/Poly}} 
\newcommand{\qsat}{\mathbf{QSAT}} 
\newcommand{\pp}{\mathbf{P}}
\newcommand{\ppp}{{\pp\pp}}

\newcommand{\nc}{\mathbf{NC}}
\newcommand{\pspa}{\mathbf{PSPACE}}

\newcommand{\eps}{{\varepsilon}}
\newcommand{\Pro}{{\mathbb{P}}}
\newcommand{\barf}{{\bar{f}}}

\newcommand{\R}{\mathbb{R}}

\newcommand{\C}{\mathbb{C}}
\newcommand{\N}{\mathbb{N}}
\newcommand{\Z}{\mathbb{Z}}
\newcommand{\bO}{\mathbf{O}}
\newcommand{\fii}{\varphi}

\newcommand{\Zn}{\Z^n}

\newcommand{\Rn}{\R^n}

\newcommand{\Cn}{\C^n}

\newcommand{\Cs}{\C^*}
\newcommand{\Rs}{\R^*}
\newcommand{\cA}{{\mathcal{A}}}
\newcommand{\hA}{{\hat{\cA}}}
\newcommand{\cB}{{\mathcal{B}}}
\newcommand{\cC}{{\mathcal{C}}}

\newcommand{\cF}{{\mathcal{F}}}
\newcommand{\cG}{{\mathcal{G}}}

\newcommand{\cV}{{\mathcal{V}}}
\newcommand{\re}{{\mathbf{Re}}}

\renewcommand{\qed}{$\blacksquare$}
\newcommand{\dia}{$\diamond$}

\newcommand{\size}{\mathrm{size}}
\newcommand{\conv}{\mathrm{Conv}}
\newcommand{\aff}{\mathrm{Aff}}

\newcommand{\supp}{\mathrm{Supp}}

\newcommand{\sign}{\mathrm{sign}}

\begin{document}
\conferenceinfo{SNC'09,} {August 3--5, 2009, Kyoto, Japan.}
\CopyrightYear{2009}
\crdata{978-1-60558-664-9/09/08}

\title{Optimization and $\text{\scalebox{1.8}[2.2]{$\np_\R$}}$-Completeness 
\\of Certain Fewnomials}
\numberofauthors{3}
\author{
\alignauthor
Philippe P\'ebay\titlenote{Supported by the United States Department of Energy, Office of Defense
Programs.}\\
\affaddr{Sandia National Laboratories}\\
\affaddr{PO Box 969, MS 9159}\\
\affaddr{Livermore, CA\ 94551}\\
\affaddr{USA}\\
\email{pppebay@sandia.gov}\\
\alignauthor
J.\ Maurice Rojas\titlenote{
Partially supported by NSF CAREER grant
DMS-0349309, the United States Department of Energy (Office of
Defense Programs), MSRI, and the Wenner Gren Foundation.}\\
\affaddr{TAMU 3368}\\
\affaddr{Department of Mathematics}\\
\affaddr{Texas A\&M University}\\
\affaddr{College Station, Texas \ 77843-3368}\\
\affaddr{USA}\\
\email{rojas@math.tamu.edu}
\alignauthor
David C.\ Thompson\titlenote{Supported by the United States Department of Energy, Office of Defense
Programs.\\[2pt]
Sandia is a multiprogram laboratory operated by Sandia
Corporation, a Lockheed Martin Company, for the United States
Department of Energy under contract DE-AC04-94AL85000.}
\affaddr{Sandia National Laboratories}\\
\affaddr{PO Box 969, MS 9159}\\
\affaddr{Livermore, CA\ 94551}\\
\affaddr{USA}\\
\email{dcthomp@sandia.gov}
}

\date{\today}

\maketitle

\noindent
{\scriptsize We dedicate this paper to Tien-Yien Li on the occasion of his 
65th birthday. Happy 65, TY!.}

\noindent
\begin{abstract}
We give a high precision polynomial-time approximation scheme for the 
supremum of any honest $n$-variate $(n+2)$-nomial with a constant term, 
allowing real exponents as well as real coefficients. Our complexity bounds 
count field operations and inequality checks, and are polynomial in $n$ and 
the {\bf logarithm} of a certain condition number. For the special case of 
polynomials (i.e., integer exponents), the log of our condition number is 
sub-quadratic in the sparse size. The best previous complexity bounds were 
exponential in the sparse size, even for $n$ fixed. Along the way, we 
partially extend the theory of $\cA$-discriminants to real exponents and 
exponential sums, and find new and natural $\np_\R$-complete problems. 
\end{abstract} 

\vspace{1mm} 
\noindent 
\scalebox{.98}[1]{{\bf Categories and Subject Descriptors:} 
F.2.1 
[Analysis of}\linebreak 
\scalebox{.9}[1]{Algorithms and Problem Complexity]: Numerical Algorithms 
and}\linebreak 
Problems---{\em Computations on 
polynomials}; G.1.5 {[Numerical Analysis]} Roots of 
Nonlinear Equations---{\em Iterative methods} 

\noindent 
{\bf General Terms:} Algorithms, Performance, Theory 


\noindent 
\section{Introduction and Main Results} 

\noindent
Maximizing or minimizing polynomial functions is a central problem in 
science and engineering. 
Typically, the polynomials have an 
underlying structure, e.g., sparsity, small expansion with respect to a 
particular basis, invariance with respect to a group action, etc. 
In the setting of sparsity, Fewnomial Theory \cite{kho} has succeeded 
in establishing bounds for the number of real solutions (or real 
extrema) that depend just on the number of monomial terms. 
However, the current general complexity bounds for real solving and 
nonlinear optimization (see, e.g., \cite{bpr,eldin,parrilo}) are 
still stated in terms of degree and number of variables, and all but ignore 
any finer input structure. In this paper, we present new speed-ups for the 
optimization of certain sparse multivariate polynomials, extended to allow 
real exponents as well. Along the way, we also present two 
new families of problems that are $\np_\R$-complete, i.e., the analogue of 
$\np$-complete for the {\bf BSS model over $\R$}. (The BSS model,  
derived in the 1980s by Blum, Shub, and Smale \cite{bss}, is a 
generalization of the classical Turing model of computation with an eye toward 
unifying bit complexity and algebraic complexity.)   

Our framework has both symbolic and numerical aspects in that 
(a) we deal with real number inputs and (b) our algorithms give either 
yes or no answers that are always correct, or numerically approximate answers 
whose precision can be efficiently tuned. Linear Programming (LP) forms 
an interesting parallel to the complexity issues we encounter. 
In particular, while LP admits polynomial-time algorithms relative to the 
Turing model, polynomial-time algorithms for linear programming relative to 
the BSS model over $\R$ (a.k.a.\ strongly polynomial-time algorithms or 
polynomial arithmetic complexity) remain unknown. Furthermore, the 
arithmetic complexity of LP appears to be linked with a fundamental 
invariant measuring the intrinsic complexity of numerical solutions:   
the {\bf condition number} (see, e.g., 
\cite{vavasisye,cucker}). Our results reveal a class of non-linear 
problems where similar subtleties arise when comparing discrete and 
continuous complexity. 

To state our results, let us first clarify some basic notation 
concerning sparse polynomials and complexity classes over $\R$. Recall that 
$\lfloor x \rfloor$ is the greatest integer not exceeding a real number $x$,  
and that $R^*$ is the multiplicative group of nonzero elements in any ring $R$. 
\newpage 
\begin{dfn} 
When $a_j\!\in\!\Rn$, the notations $a_j\!=$\linebreak
$(a_{1,j},\ldots,a_{n,j})$, 
$x^{a_j}\!=\!x^{a_{1,j}}_1\cdots x^{a_{n,j}}_n$, and $x\!=\!(x_1,\ldots,x_n)$ 
will be understood. If $f(x)\!:=\!\sum^m_{j=1} c_ix^{a_j}$ 
where $c_j\!\in\!\Rs$ for all $j$, 
and the $a_j$ are pair-wise distinct, then we call $f$ a 
{\bf (real) $\pmb{n}$-variate $\pmb{m}$-nomial}, and we define 
$\supp(f)\!:=\!\{a_1,\ldots,a_m\}$ to be the {\bf support} of $f$. 
We also let $\cF_{n,m}$ denote the set of all $n$-variate 
$\lfloor m\rfloor$-nomials\footnote{Here we allow real coefficients, 
unlike \cite{finally} where the same notation included a restriction to integer 
coefficients.} and, for any 
$m\!\geq\!n+1$, we let $\cF^*_{n,m}\!\subseteq\!\cF_{n,m}$ 
denote the subset consisting of those $f$ with $\supp(f)$ {\bf not} 
contained in any $(n-1)$-flat. We also call any 
$f\!\in\!\cF^*_{n,m}$ an {\bf honest 
$\pmb{n}$-variate $\pmb{m}$-nomial} (or {\bf honestly $\pmb{n}$-variate}). \dia
\end{dfn}
For example, the dishonestly $4$-variate trinomial\\
\mbox{}\hfill $-1+\sqrt{7}x^2_1x_2x^7_3x^3_4-e^{43} 
x^{198e^2}_1x^{99e^2}_2x^{693e^2}_3x^{297e^2}_4$ \hfill\mbox{}\\
(with support contained in a line segment) has the same supremum 
over $\R^4_+$ as the {\bf honest uni}variate trinomial\\ 
\mbox{}\hfill 
$-1+\sqrt{7}y_1-e^{43}y^{99e^2}_1$ 
\hfill\mbox{}\\ 
has over $\R_+$.
More generally, via a monomial change of variables, 
it will be natural to restrict to $\cF^*_{n,n+k}$ (with $k\!\geq\!1$) to 
study the role of sparsity in algorithmic complexity over the real numbers. 

We will work with some well-known complexity classes from the 
BSS model over $\R$ (treated fully in \cite{bcss}), so we will only briefly 
review a few definitions, focusing on a particular extension we need. 
Our underlying notion of\linebreak input size, including a 
variant of the 
condition number, is clarified in Definition \ref{dfn:cond} of Section 
\ref{sub:input} below, and\linebreak illustrated 
in Example \ref{ex:input} immediately following our first main theorem.

So for now let us just recall the following basic\linebreak 
inclusions of complexity 
classes: $\nc^1_\R\!\subsetneqq\!\pp_\R\!\subseteq\!\np_\R$ \cite[Ch.\ 19, 
Cor.\ 1, pg.\ 364]{bcss}. (The properness of the latter\linebreak 
inclusion remains a 
famous open problem, akin to the more famous classical 
$\pp\text{\scalebox{1}[.7]{$\stackrel{?}{=}$}}\np$ 
question.) Let us also recall that $\nc^k_\R$ is the family of real valued 
functions (with real inputs) 
computable by arithmetic 
circuits\footnote{This is one of $2$ times we will mention circuits
in the sense of complexity theory: Everywhere else in this paper,
our circuits will be {\bf combinatorial} objects as in Definition
\ref{dfn:ckt} below.} with size polynomial in the input size 
and depth $O\!\left(\log^k\!\left(\text{Input Size}\right)\right)$ 
(see \cite[Ch.\ 18]{bcss} for further discussion). 

To characterize a natural class of problems with efficiently 
computable numerical answers, we will define the notion of a 
{\bf High Precision Polynomial Time Approximation Scheme}: We 
let $\hpt$ denote the class 
of functions $\phi : \R^\infty\longrightarrow \R\cup\{+\infty\}$ 
such that, for any $\eps\!>\!0$, there is an algorithm guaranteed to 
approximate $\phi(x)$ to within a $1+\eps$ factor, using 
a number of arithmetic operations \linebreak 
polynomial in $\size(x)$ 
{\bf and} $\log\log\frac{1}{\eps}$.\footnote{When $\phi(x)\!=\!0$ we 
will instead require an {\bf additive} error of $\eps$ or less. 
When $\phi(x)\!=\!+\infty$ we will require the approximation to 
be $+\infty$, regardless of $\eps$.} 
Our notation is\linebreak 
inspired by the more familiar classical family of 
problems $\fpt$ (i.e., those problems admitting a {\bf Fully\linebreak 
Polynomial Time 
Approximation Scheme}), where\linebreak instead the input is 
Boolean and the complexity 
need only be polynomial in $\frac{1}{\eps}$. The complexity class $\fpt$ was 
\linebreak 
formulated in \cite{acg} and a highly-nontrivial\linebreak 
example of a problem 
admitting a $\fpt$ is counting\linebreak 
matchings in bounded degree graphs \cite{count}. 

\newpage 

\begin{rem} 
For a vector function $\phi=(\phi_1,\ldots,\phi_k) : \R^\infty \longrightarrow 
(\R\cup\infty)^k$ it will be natural to say that $\phi$ admits an 
$\hhpt$ iff each coordinate of $\phi_i$ admits an $\hhpt$. \dia 
\end{rem} 

\subsection{Sparse Real Optimization} 

\noindent 
The main computational problems we address are the\linebreak 
following. 
\begin{dfn}
Let $\R_+$ denote the positive real numbers, and let $\opt$ denote the 
problem of deciding, for a given $(f,\lambda)\!\in\!
\left(\bigcup\limits_{n\in\N}\R[x^a \; | \; a\!\in\!\Rn]\right)\times 
\R$, whether $\sup_{x\in\Rn_+} f\!\geq\!\lambda$ or 
not. Also, for any subfamily $\cF\subseteq\bigcup_{n\in\N} 
\R[x^a\; | \; a\!\in\!\Rn]$, we let
$\opt(\cF)$ denote the natural restriction of $\opt$ 
to inputs in $\cF$. Finally, we let $\oopt$ (resp.\ $\oopt(\cF)$)  
denote the obvious functional analogue of $\opt$ (resp.\ 
$\opt(\cF)$) where (a) the input is instead $(f,\eps)\!\in\!
\left(\bigcup\limits_{n\in\N}\R[x^a \; | \; a\!\in\!\Rn]\right)\times 
\R_+$ and (b) the output is instead a pair\\ 
\mbox{}\hfill $(\bar{x},\bar{\lambda})\!\in\!(\R_+\cup\{0,+\infty\})^n\times 
(\R\cup\{+\infty\})$\hfill \mbox{}\\ 
with $\bar{x}\!=\!(\bar{x}_1,\ldots,\bar{x}_n)$ (resp.\ 
$\bar{\lambda}$) an $\hhpt$ for $x^*$ (resp.\ $\lambda^{*}$) 
where $\lambda^*\!:=\!\sup_{x\in\Rn_+} f\!=\!\lim_{x\rightarrow x^*}f(x)$ 
for some $x^{*}\!=\!(x^*_1,\ldots,x^*_n)\!\in\!(\R_+\cup\{0,+\infty\})^n$. 
\end{dfn}
\begin{rem} 
Taking logarithms, it is clear that our\linebreak problems above   
are equivalent to maximizing a function of the form $g(y)\!=\!
\sum^m_{i=1}c_ie^{a_i\cdot y}$ over $\Rn$. When convenient, we will use 
the latter notation but, to draw parallels with 
the algebraic case, we will usually speak of ``polynomials'' with real 
exponents. \dia 
\end{rem} 

We will need to make one final restriction when\linebreak optimizing 
$n$-variate $m$-nomials: we let $\cF^{**}_{n,n+k}$ denote the 
subset of $\cF^*_{n,n+k}$ consisting of those $f$ with $\supp(f)\!\ni\!
\bO$. While technically convenient, this restriction is also natural 
in that level sets of $(n+k)$-nomials in $\cF^{**}_{n,n+k}$ become  
zero sets of $(n+k')$-nomials with $k'\!\leq\!k$. 

We observe that checking whether the zero set of an $f\!\in\!\R[x_1,
\ldots,x_n]$ is nonempty (a.k.a.\ the {\bf real (algebraic)\linebreak 
feasibility 
problem}) is equivalent to checking whether the maximum of 
$-f^2$ is $0$ or greater.  So it can be argued 
that the $\np$-hardness (and $\np_\R$-hardness) of $\opt$ has 
been known at least since the 1990s \cite{bcss}. However, it appears that no 
sharper complexity upper bounds in terms of sparsity were known earlier. 
\begin{thm}
\label{THM:BIG} 
We can efficiently optimize $n$-variate\linebreak $(n+k)$-nomials 
over $\Rn_+$ for $k\!\leq\!2$. Also, for $k$ a slowly growing function 
of $n$, optimizing $n$-variate $(n+k)$-nomials over $\Rn_+$ is 
$\np$-hard. More precisely:
\begin{enumerate} 
\addtocounter{enumi}{-1}
\item{Both $\opt\!\left(\bigcup_{n\in\N} \cF^{**}_{n,n+1}\right)$ and 
$\oopt\!\left(\bigcup_{n\in\N} \cF^{**}_{n,n+1}\right)$ 
are in $\nc^1_\R$. } 
\item{$\opt\!\left(\bigcup_{n\in\N}\cF^{**}_{n,n+2}\right)\!\in\!\pp_\R$  
and $\oopt\!\left(\bigcup_{n\in\N}\cF^{**}_{n,n+2}\right)\!\in\!\hpt$. } 
\item{For any fixed $\delta\!>\!0$, \scalebox{.85}[1]
{$\opt\text{\raisebox{-.3cm}{\scalebox{1}[2.3]{$($}}}
\bigcup\limits_{\substack{n\in\N\\ 
0<\delta'<\delta}} \cF^{**}_{n,n+n^{\delta'}}
\cap\R[x_1,\ldots,x_n]\text{\raisebox{-.3cm}{\scalebox{1}[2.3]{$)$}}}$} 

\vspace{-.6cm}
is $\np_\R$-complete.} 
\end{enumerate} 
\end{thm} 

\begin{ex}
\label{ex:input} 
Suppose $\eps\!>\!0$. 
A very special case of \linebreak 
Assertion (1) of Theorem \ref{THM:BIG} then implies 
that we can\linebreak approximate within a factor of $1+\eps$ --- 
for any real nonzero $c_1,\ldots,c_{n+2}$ and  
$D$ ---  the maximum of the function $f(x)$ defined to be\\ 
\scalebox{.8}[1]{$c_1+c_2(x^D_1\cdots x^{D^n}_n) + 
c_3 x^{2D}_1\cdots x^{2^n D^n}_n +\cdots+c_{n+2}x^{(n+1)D}_1\cdots 
x^{(n+1)^nD^n}_n$,}\\  
using a number of arithmetic operations linear in\\
\mbox{}\hfill $n^2(\log(n)+\log D)+\log\log\frac{1}{\eps}$. 
\hfill\mbox{}\\ 
The best previous results in the algebraic setting 
(e.g., the critical points method as detailed in \cite{eldin}, or 
by combining \cite{bpr} and the efficient numerical approximation results of 
\cite{mp98}) would yield a bound polynomial in\\
\mbox{}\hfill $n^nD^{n}+\log\log\frac{1}{\eps}$ \hfill \mbox{}\\ 
instead, and only under the assumption that $D\!\in\!\N$. 
Alternative approaches via semidefinite programming also appear to result in 
complexity bounds superlinear in $n^nD^{n}$ (see, e.g., 
\cite{parrilo,lassere,niesparse,kojima}), and still require 
$D\!\in\!\N$. Moving to Pfaffian/Noetherian function 
techniques, \cite{gv} allows arbitrary real $D$ but still yields an arithmetic 
complexity bound exponential in $n$. It should of course be pointed out that 
the results of \cite{bpr,mp98,eldin,parrilo,lassere,niesparse,kojima,gv} apply 
to real polynomials in complete generality. \dia  
\end{ex}
We thus obtain a significant speed-up for a particular class of 
analytic functions, laying some preliminary groundwork for improved 
optimization of $(n+k)$-nomials with $k$\linebreak 
arbitrary. Our advance is possible 
because, unlike earlier methods which essentially revolved around commutative 
\linebreak 
algebra (and were more suited to complex algebraic\linebreak geometry), we are 
addressing a real analytic problem with real analytic tools. 
Theorem \ref{THM:BIG} is proved in Section \ref{sub:thresh}\linebreak 
below.
Our main new technique, which may be of\linebreak independent interest, is 
an extension of $\cA$-discriminants\linebreak (a.k.a.\ sparse 
discriminants) to real 
exponents (Theorem \ref{THM:DISC} of Section \ref{sub:disc}). 

Our algorithms are quite implementable (see Algorithm \ref{algor:fsup} 
of Section \ref{sub:thresh}) and derived via a 
combination of\linebreak tropical geometric ideas and $\cA$-discriminant 
theory, both extended to real exponents. In particular, for $n$-variate 
$(n+1)$-nomials, a simple change of variables essentially tells us that 
tropical geometry rules (in the form of {\bf Viro diagrams} \cite[Ch.\ 5,
pp.\ 378--393]{gkz94}, but extended to real exponents), and in the case at 
hand this means that one can 
compute extrema by checking inequalities involving the coefficients (and 
possibly an input $\lambda$). 
Tropical geometry still applies to the $n$-variate $(n+2)$-nomial case, but 
only after one evaluates the sign of a particular generalized 
$\cA$-discriminant.\footnote{For $n$-variate $(n+3)$-nomials, knowing the sign 
of a discriminant is no longer sufficient, and efficient optimization 
still remains an open problem. Some of the intricacies are detailed in 
\cite{drrs,reu08}.} More precisely, an $n$-variate $m$-nomial $f$ (considered 
as a function on $\Rn_+$) with bounded supremum $\lambda^*$ must 
attain the value $\lambda^*$ at a critical point of $f$ in the nonnegative 
orthant. In particular, the nonnegative zero set of $f-\lambda^*$ must be 
degenerate, and thus we can attempt to solve for $\lambda^*$ (and a 
corresponding maximizer) if we have a sufficiently tractable notion of 
discriminant to work with. 

So our hardest case reduces to (a) finding efficient 
formulas for discriminants of $n$-variate $m$-nomials and\linebreak (b) 
efficiently 
detecting unboundedness for $n$-variate $m$-\linebreak 
nomials. When $m\!=\!n+2$, (a) 
fortuitously admits a solution, based on a nascent theory developed further 
in \cite{evy}. We can also reduce Problem (b) to Problem (a) via 
some\linebreak tropical geometric tricks. So our development ultimately hinges 
deriving an efficient analogue of discriminant polynomials for discriminant 
varieties that are no longer algebraic. 
\newpage 

\begin{ex} 
Consider the trivariate pentanomial 
$f:=$\linebreak
$c_1+c_2x^{999}_1+c_3x^{73}_1x^{\sqrt{363}}_3
+c_4x^{2009}_2+c_5x^{74}_1x^{108e}_2x_3$, 
with $c_1,\ldots,$\linebreak 
$c_4\!<\!0$ and $c_5\!>\!0$. Theorem \ref{thm:bigger} 
of Section \ref{sub:degen} then\linebreak easily implies that 
$f$ attains a maximum of $\lambda^*$ on $\R^3_+$  
iff $f-\lambda^*$ has a degenerate root in $\R^3_+$. Via  
Theorem \ref{THM:DISC} of Section \ref{sub:disc} below,  
the latter occurs iff\\
\mbox{}\hfill $b^{b_5}_5(c_1-\lambda^*)^{b_1} c^{b_2}_2
c^{b_3}_3 c^{b_4}_4-b^{b_1}_1 b^{b_2}_2 b^{b_3}_3 b^{b_4}_4 c^{b_5}_5$
\hfill\mbox{}\\  
vanishes, where $b\!:=\!(b_1,b_2,b_3,b_4,-b_5)$ is any generator of the 
kernel of the map $\varphi : \R^5\longrightarrow \R^4$ defined by 
the matrix\\
\mbox{}\hfill \scalebox{.7}[.7]{$\begin{bmatrix}1 & 1 & 1 & 1 & 1\\
0 & 999        & 73 & 0    & 74\\
0 & 0          & 0  & 2009 & 108e\\
0 & \sqrt{363} & 0  & 0    & 1 \end{bmatrix}$,}\hfill\mbox{}\\  
normalized so that $b_5\!>\!0$. 
In particular, such a $b$ can be computed easily via $5$  
determinants of $4\times 4$ submatrices (via Cramer's Rule), and we thus see 
that $\lambda^*$ is nothing more than $c_1$ minus a monomial 
(involving real exponents) in $c_2,\ldots,c_5$. Via the 
now classical fast algorithms for approximating $\log$ and 
$\exp$ \cite{brent}, real powers of real numbers (and thus $\lambda^*$) 
can be efficiently approximated. Similarly, deciding whether $\lambda^*$  
exceeds a given $\lambda$ reduces to checking an inequality involving real 
powers of positive numbers. \dia
\end{ex}

\subsection{Related Work} 

\noindent 
The computational complexity of numerical analysis continues to be 
an active area of research, both in theory and in practice. On the 
theoretical side, the BSS model over $\R$ has proven quite useful 
for setting a rigourous foundation. While this model involves exact 
arithmetic and field operations, there are many results building upon this 
model that elegantly capture round-off error and numerical 
conditioning (see, e.g., \cite{cuckersmale,burgisser}). 
Furthermore, results on $\pp_\R$ and $\np_\R$ do ultimately impact classical 
complexity classes. For instance, the respective {\bf Boolean parts} of these 
complexity classes, $\bp(\pp_\R)$ and $\bp(\np_\R)$, are defined as 
the respective restrictions of $\pp_\R$ and $\np_\R$ to integer inputs. 
While the best known bounds for these Boolean parts are still rather loose 
--- \\ 
\mbox{}\hfill $\ppoly\!\subseteq\!\bp(\pp_\R)\!\subseteq\!\pspoly$ 
\cite{burgisser}, 
\hfill\mbox{}\\
\mbox{}\hfill $\npoly\!\subseteq\!\bp(\np_\R)\!\subseteq\!\ch$ 
\cite{burgisser},\hfill\mbox{}\\ 
--- good algorithms for the BSS model and good algorithms for the 
Turing model frequently inspire one another, e.g., \cite{realkoiran,bpr}. 
We recall that $\ppoly$, referred to as {\bf non-uniform 
polynomial-time}, consists of those decision problems solvable by a 
non-uniform family of circuits\footnote{i.e., there is no restriction 
on the power of the algorithm specifying the circuit 
for a given input size} of size polynomial in the input.  
$\ch$ is the {\bf counting hierarchy} $\ppp\cup \ppp^{\ppp}
\cup \ppp^{\ppp^\ppp}\cup \cdots$, which happens to be contained 
in $\pspa$ (see \cite{burgisser} and the references therein). 

Let us also point out that the number of natural\linebreak 
problems known to be 
$\np_\R$-complete remains much smaller than the number of natural 
problems known to be $\np$-complete: deciding the 
existence of a real roots for\linebreak multivariate polynomials (and 
various subcases involving\linebreak 
quadratic systems or single quartic 
polynomials) \cite[Ch.\ 5]{bcss}, linear programming feasibility \cite[Ch.\ 
5]{bcss}, and bounding the real dimension of algebraic sets \cite{realkoiran} 
are the main representative $\np_\R$-complete problems.\linebreak  
Optimizing $n$-variate $(n+n^\delta)$-nomials (with $\delta\!>\!0$ 
fixed and $n$ unbounded), and the corresponding feasibility 
problem (cf.\ Corollary \ref{COR:RFEAS} below), now join this 
short list. 

While sparsity has been profitably explored in 
the context of interpolation (see, e.g., \cite{ky,gll}) and factorization over 
number fields \cite{lenstra,kk,aks}, 
it has been mostly ignored in numerical analysis (for nonlinear polynomials) 
and the study of the BSS model over $\C$ and $\R$. 
For instance, there appear to be no earlier 
published complexity upper bounds of the form 
$\opt\left(\cF_{1,m}\right)\!\in\!\pp_\R$ (relative to the sparse encoding) 
for any $m\!\geq\!3$, in spite of beautiful recent work in 
semi-definite programming (see, e.g., \cite{lassere,niesparse,kojima}) 
that begins to address the optimization of sparse multivariate polynomials
over the real numbers. In particular, while the latter papers give significant 
practical speed-ups over older techniques such as resultants and Gr\"obner 
bases, the published complexity bounds are still exponential (relative to the 
sparse encoding) for $n$-variate $(n+2)$-nomials, and require the assumption 
of integer exponents. 

We can at least obtain a glimpse of sparse optimization beyond $n$-variate 
$(n+2)$-nomials by combining our framework with an earlier result from 
\cite{rojasye}. The proof is in Section \ref{sub:4}.  
\begin{cor} 
\label{COR:4} \mbox{}\\
(0) Using the same notion of input size as for $\oopt$ (cf.\linebreak 
\mbox{}\hspace{.6cm}Definition \ref{dfn:cond} below), the positive roots of 
any real\linebreak 
\mbox{}\hspace{.6cm}trinomial in $\cF_{1,3}\cap\R[x_1]$ admit an $\hhpt$.\\
(1) $\opt(\cF^{**}_{1,4}\cap\R[x_1])\!\in\!\pp_\R$ and 
$\oopt(\cF^{**}_{1,4}\cap\R[x_1])\!\in$\linebreak
\mbox{}\hspace{.6cm}$\hpt$. 
\end{cor} 

As for earlier complexity lower bounds for $\opt$ in terms of sparsity,
we are unaware of any. For instance, it is not even known whether
$\opt(\R[x_1,\ldots,x_n])$ is $\np_\R$-hard for some fixed $n$ 
(relative to the sparse encoding). 

The paper \cite{finally}, which deals exclusively with decision 
problems (i.e., yes/no answers) and bit complexity (as opposed to 
arithmetic complexity), is an important precursor to the present work. 
Here, we thus expand the context to real coefficient and 
real exponents, work in the distinct setting of optimization, and derive 
(and make critical use of) a new tool: 
generalized $\cA$-discriminants for exponential sums. As a consequence, 
we are also able to extend some of the complexity lower bounds from 
\cite{finally} as follows. (See Section \ref{sub:thresh} for the proof.) 
\begin{dfn}  
Let $\feas_\R$ (resp.\ $\feas_+$) denote the problem of deciding whether an 
arbitrary system of equations from 
$\bigcup_{n\in\N} \R[x^a\; | \; a\!\in\!\Rn]$ has a real root (resp.\ 
root with all coordinates positive). Also, for any collection $\cF$ of tuples 
chosen from $\bigcup_{k,n\in\N}(\R[x^a\; | \; a\!\in\!\Rn])^k$, we let 
$\feas_\R(\cF)$ (resp.\ $\feas_+(\cF)$) denote the natural restriction of 
$\feas_\R$ (resp.\ $\feas_+$) to inputs in $\cF$. \dia 
\end{dfn} 
\begin{cor} 
\label{COR:RFEAS}  
For any $\delta\!>\!0$,\\ 
\mbox{}\hfill 
$\feas_\R\text{\raisebox{-.3cm}{\scalebox{1}[2.3]{$($}}}
\bigcup\limits_{\substack{n\in\N\\ 0<\delta'<\delta}}
\cF^{**}_{n,n+n^{\delta'}}\cap\R[x_1,\ldots,x_n] 
\text{\raisebox{-.3cm}{\scalebox{1}[2.3]{$)$}}} 
$\hfill \mbox{}\\ 
and \\
\mbox{}\hfill 
$\feas_+
\text{\raisebox{-.3cm}{\scalebox{1}[2.3]{$($}}} 
\bigcup\limits_{\substack{n\in\N\\ 0<\delta'<\delta}}
\cF^{**}_{n,n+n^{\delta'}}\cap\R[x_1,\ldots,x_n]
\text{\raisebox{-.3cm}{\scalebox{1}[2.3]{$)$}}} 
$\hfill \mbox{}\\ 
are each $\np_\R$-complete.  
\end{cor} 

\section{Background}  
\label{sec:back} 
\subsection{Input Size} 
\label{sub:input} 
To measure the complexity of our algorithms,  
let us fix the following definitions for input size and condition number. 
\begin{dfn} 
\label{dfn:cond} 
Given any subset $\cA\!=\!\{a_1,\ldots,a_m\}\!\subset\!\Rn$ of cardinality 
$m$, let us define $\hA$ to be the $(n+1)\times m$ matrix whose $j\thth$ 
column is $\{1\}\times a_j$, and $\beta_J$ the absolute value of the 
determinant of the submatrix of $\hA$ consisting of those columns 
of $\hA$ with index in a subset $J\!\subseteq\!\{1,\ldots,m\}$ 
of cardinality $n+1$. Then,  
given any $f\!\in\!\cF^*_{n,m}$ written $f(x)\!=\!\sum^{m}_{i=1}
c_ix^{a_i}$, we define its {\bf condition number}, $\cC(f)$, to be\\ 
\mbox{}\hfill $\left(\prod\limits^{m}_{i=1}\max\!\left
\{3,|c_i|,\frac{1}{|c_i|}\right\}\right)\times 
\prod\limits_{\substack{J\subseteq\{1,\ldots,m\}\\ \#J=n+1}} 
\max^*\!\left(3,|\beta_J|,\frac{1}{|\beta_J|}\right)$,\hfill\mbox{}\\ 
where $\max^*(a,b,c)$ is $\max\{a,b,c\}$ or $a$, according as 
$\max\{b,c\}$ is finite or not. 

Throughout this paper, we will use the following notions of input size for 
$\opt$ and $\oopt$: 
The size of any\linebreak 
instance $(f,\lambda)$ of $\opt$ (resp.\ an instance 
$(f,\eps)$ of $\oopt$) is $\log\!\left(\max^*\left(3,|\lambda|,\frac{1}
{|\lambda|}\right)\right)+\log \cC(f)$
(resp.\ $\log \cC(f)$). \dia  
\end{dfn} 

While our definition of condition number may appear\linebreak unusual, it is 
meant to concisely arrive at two important properties: 
(1) $\log \cC(f)$ is polynomial in $n\log \deg f$ when 
$f\!\in\!\cF_{n,n+k}\cap\R[x_1,\ldots,x_n]$ and $k$ is fixed, (2) $\cC(f)$ 
is closely\linebreak 
related to an underlying discriminant (see Theorem \ref{THM:DISC}\linebreak 
below) that dictates how much numerical accuracy we will\linebreak 
\scalebox{.89}[1]{need to solve 
$\oopt$. We also point out that for $f\!\in\!\Z[x_1,\ldots,x_n]$,}\linebreak  
it is easy to show that $\log \cC(f)\!=\!O(nS(f))$ where $S(f)$ is the 
{\bf sparse size} of $f$, i.e., $S(f)$ is the number of bits needed to write 
down the monomial term expansion of $f$. For\linebreak sufficiently sparse 
polynomials, 
algorithms with\linebreak 
complexity polynomial in $S(f)$ are much faster than those with 
complexity polynomial in $n$ and $\deg(f)$. \cite{lenstra,kk,aks,ky,gll,
finally} provide other interesting\linebreak  
examples of algorithms with complexity polynomial in $S(f)$. 

\subsection{Tricks with Exponents} 
\label{sub:exp} 

\noindent 
A simple and useful change of variables is to use\linebreak 
monomials in new variables. 
\begin{dfn}
For any ring $R$, let $R^{m\times n}$ denote the set of $m\times n$ matrices
with entries in $R$. For any $M\!=\![m_{ij}]\!\in\!\R^{n\times n}$ 
and $y\!=\!(y_1,\ldots,y_n)$, we define the formal expression 
$y^M\!:=\!(y^{m_{1,1}}_1\cdots y^{m_{n,1}}_n,\ldots, 
y^{m_{1,n}}_1\cdots y^{m_{n,n}}_n)$. We call the substitution 
$x\!:=\!y^M$ a {\bf monomial change of variables}. Also, for 
any $z\!:=\!(z_1,\ldots,z_n)$, we let $xz\!:=\!(x_1z_1,\ldots,x_nz_n)$. 
\linebreak 
Finally, let $\gln(\R)$ denote the group of all 
invertible matrices in $\R^{n\times n}$. \dia 
\end{dfn} 

\begin{prop} 
\label{prop:monochange} 
(See, e.g., \cite[Prop.\ 2]{tri}.) 
For any $U,V\!\in\!\R^{n\times n}$, we have the formal identity\\ 
\mbox{}\hfill $(xy)^{UV}\!=\!(x^U)^V(y^U)^V$.\hfill\mbox{}\\ 
Also, if $\det U\!\neq\!0$, then the function 
$e_U(x)\!:=\!x^U$ is an\linebreak analytic 
automorphism of $\Rn_+$, and preserves smooth points and singular 
points of 
positive zero sets of analytic functions.
Finally, $U\!\in\!\gln(\R)$ 
implies that $e^{-1}_U(\Rn_+)\!=\!\Rn_+$ and that $e_U$ 
maps distinct open orthants of $\Rn$ to distinct open orthants of $\Rn$. \qed 
\end{prop} 

A consequence follows: Recall that the {\bf affine span}
of a point set $\cA\!\subset\!\Rn$, $\aff \cA$, is the set of real
linear\linebreak combinations $\sum_{a\in\cA} c_a a$ satisfying 
$\sum_{a\in\cA}c_a\!=\!0$. To\linebreak 
optimize an $f\!\in\!
\cF^{**}_{n,n+1}$ it will help to have a much simpler canonical form.
In what follows, we use $\#$ for set cardinality and 
$e_i$ for the $i\thth$ standard basis 
vector of $\Rn$.  
\vfill\eject 

\begin{cor}
\label{cor:can}
For any $f\!\in\!\cF^{**}_{n,n+1}$ we can compute 
$c\!\in\!\R$ and $\ell\!\in\!\{0,\ldots,n\}$ within $\nc^1_\R$ such that\\ 
\mbox{}\hfill $\barf(x)\!:=\!c+x_1+\cdots+x_\ell-x_{\ell+1}-\cdots -x_n$ 
\hfill\mbox{}\\
satisfies:\\ 
\mbox{}\hspace{.4cm}(1) $\bar{f}$ and $f$ have exactly the same number of 
positive\\ 
\mbox{}\hspace{1cm}coefficients, and\\ 
\mbox{}\hspace{.4cm}(2) 
$\bar{f}\!\left(\Rn_+\right)\!=\!f\!\left(\Rn_+\right)$. 
\end{cor}

\noindent 
{\bf Proof:} Suppose $f$ has support $\cA\!=\!\{0,a_2,\ldots,a_{n+1}\}$ 
and corresponding coefficients $c_1,\ldots,c_{n+2}$. 
Letting $B$ denote the $n\times n$ matrix whose $i\thth$ column is $a_{i+1}$, 
Proposition \ref{prop:monochange}, via the substitution 
$x\!=\!y^{B^{-1}}$, tells us that we may assume that 
$f$ is of the form $c_1+c_2x_1+\cdots+c_{n+1}x_n$. Moreover, 
to obtain $\bar{f}$, we need only perform a suitable positive rescaling and 
reordering of the variables. In summary, $c$ is simply the constant 
term of $f$ and $\ell$ is the number of positive 
coefficients not belonging to the constant term --- both of 
which can be computed simply by a search and a sort clearly belonging to 
$\nc^1_\R$. \qed 

\smallskip 
\noindent 
Note that we don't actually need to compute $B^{-1}$ to obtain $\ell$: 
$B^{-1}$ is needed only for the proof of our corollary.  

A final construction we will need is the notion of a\linebreak  
{\bf generalized} Viro diagram. Recall that a {\bf triangulation} of a 
point set $\cA$ is simply a simplicial complex $\Sigma$ whose vertices lie in 
$\cA$. We say that a triangulation of $\cA$ is {\bf induced by a lifting} iff 
it its simplices are exactly the domains of linearity for some function 
that is convex, continuous, and piecewise linear on the convex hull 
of\footnote{i.e., smallest convex set containing...} 
$\cA$.  
\begin{dfn} 
Suppose $\cA\!\subset\!\Rn$ is finite, $\dim \aff \cA\!=\!n$, and $\cA$
is equipped with a triangulation $\Sigma$ induced by a lifting {\bf and} a
function $s : \cA \longrightarrow \{\pm\}$ which we will call a
{\bf distribution of signs for $\cA$}. We then define a
piece-wise linear manifold  --- the {\bf generalized Viro diagram} 
$\cV_\cA(\Sigma,s)$ --- in the following local manner: For any $n$-cell 
$C\!\in\!\Sigma$,
let $L_C$ be the convex hull of the set of midpoints of edges of
$C$ with vertices of opposite sign, and then define
$\cV_\cA(\Sigma,s)\!:=\!\bigcup\limits_{C \text{ an } n\text{-cell}}
L_C$. When $\cA\!=\!\supp(f)$ and $s$ is the corresponding sequence of
coefficient signs, then we also call $\cV_{\Sigma}(f)\!:=\!\cV_\cA(\Sigma,s)$
the {\bf (generalized) Viro diagram of $f$}. \dia
\end{dfn}

\noindent
We use the appelation ``generalized'' since, to the best of our knowledge, 
Viro diagrams have only been used in the special case $\cA\!\subset\!\Zn$ 
(see, e.g., Proposition 5.2 and Theorem 5.6 of 
\cite[Ch.\ 5, pp.\ 378--393]{gkz94}). We give examples of 
Viro diagrams in Section \ref{sub:degen} below. 

\subsection{Generalized Circuit Discriminants and\\ Efficient Approximations} 
\label{sub:disc} 

\noindent 
Our goal here is to extract an extension of $\cA$-discriminant theory 
sufficiently strong to prove our main results.
\begin{dfn} 
\label{dfn:adisc} 
Given any $\cA\!=\!\{a_1,\ldots,a_m\}\!\subset\!\Rn$ of 
cardinality $m$ and $c_1,\ldots,c_m\!\in\!\Cs$, we 
define $\nabla_\cA\!\subset\!\Pro^{m-1}_\C$ --- 
the {\bf generalized $\cA$-discriminant variety} --- 
to be the closure of the 
set of all $[c_1:\cdots :c_m]\!\in\!\Pro^{m-1}_\C$ such that 
$g(x)\!=\!\sum^m_{i=1} c_ie^{a_i\cdot y}$ has a degenerate root in $\Cn$. 
In particular, we call $f$ an {\bf $n$-variate exponential $m$-sum}. 
\dia 
\end{dfn} 

To prove our results, it will actually suffice to deal with a small 
subclass of $\cA$-discriminants.
\begin{dfn} 
\label{dfn:ckt} 
We call $\cA\!\subset\!\Rn$ a {\bf (non-degenerate) 
circuit}\footnote{This terminology comes from matroid theory and
has nothing to do with circuits from complexity theory.} 
iff $\cA$ is affinely dependent, but every proper subset of $\cA$ is affinely
independent. Also, we say that $\cA$ is a {\bf degenerate circuit} iff
$\cA$ contains a point $a$ and a proper subset $\cB$ such that $a\!\in\!\cB$, 
$\cA\setminus a$ is affinely independent, and $\cB$ is a non-degenerate 
circuit. \dia 
\end{dfn}

\noindent 
For instance, both \epsfig{file=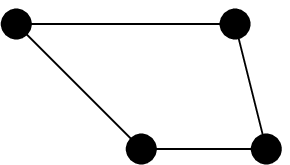,height=.35cm} and
\epsfig{file=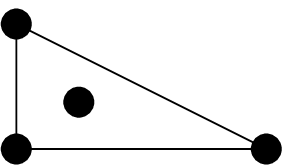,height=.35cm} are circuits, but
\epsfig{file=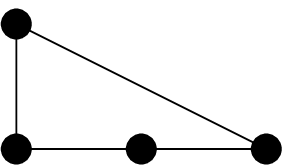,height=.35cm} is a degenerate circuit.
In general, for any degenerate circuit $\cA$, the subset $\cB$ named above is 
always unique. 

\begin{dfn} 
\label{dfn:chamber} 
For any $\cA\!\subset\!\Rn$ of cardinality $m$, let $\cG_\cA$ 
denote the set of all $n$-variate exponential $m$-sums with support $\cA$. \dia 
\end{dfn} 

There is then a surprisingly succinct description for $\nabla_\cA$ when $\cA$ 
is a non-degenerate circuit. The theorem below is inspired by 
\cite[Prop.\ 1.2, pg.\ 217]{gkz94} and \cite[Prop.\ 1.8, Pg.\ 274]{gkz94} 
--- important precursors that covered the special case of integral  
exponents. 
\begin{thm}
\label{THM:DISC} 
Suppose $\cA\!=\!\{a_1,\ldots,a_{n+2}\}\!\subset\!\Rn$ is a non-degenerate 
circuit, and let $b\!:=\!(b_1,\ldots,b_{n+2})$ where 
$b_i$ is $(-1)^i$ times the determinant 
of the matrix with columns $1\times a_1,\ldots,\widehat{1\times a_i},
\ldots,a_{n+2}$ ($\widehat{(\cdot)}$ denoting omission). Then:  
\begin{enumerate}
\item{$\nabla_\cA\!\subseteq\!\left\{[c_1:\cdots:c_{n+2}]\!\in\!\Pro^{n+1}_\C\;
: \; \prod\limits^{n+2}_{i=1} \left|\frac{c_i}{b_i}\right|^{b_i}
\!=\!1\right\}$. Also, $(b_1,\ldots,b_{n+2})$ can be computed in 
$\nc^2_\R$. }  
\item{There is a $[c_1:\cdots:c_{n+2}]\!\in\!\Pro^{n+1}_\R$ 
with\\
\mbox{}\hspace{1cm}(i) $\sign(c_1b_1)\!=\cdots=\!\sign(c_{n+2}b_{n+2})$\\
and\\    
\mbox{}\hspace{1cm}(ii) $\prod\limits^{n+2}_{i=1} 
(\sign(b_ic_i)c_i/b_i)^{\sign(b_ic_i)b_i}\!=\!1$\\ 
iff the real zero set of 
$g(y)\!:=\!\sum^{n+2}_{i=1}c_i e^{a_i\cdot y}$ 
contains a degenerate point $\zeta$. In particular, any such 
$\zeta$ satisfies $e^{a_i\cdot \zeta}\!=\!\sign(b_1c_1)b_i/c_i$ 
for all $i$, and thus the real zero set of $g$ has at most one degenerate 
point.} 
\end{enumerate}
\end{thm}

\noindent 
Theorem \ref{THM:DISC} is proved in Section \ref{sec:proofs} below. 

We will also need a variant of a family of fast algorithms discovered 
independently by Brent and Salamin.
\begin{brent} 
\cite{brent,salamin} 
Given any\linebreak positive $x,\eps\!>\!0$, we can approximate 
$\log x$ and $\exp(x)$ within a factor of $1+\eps$ using just 
$O\!\left(|\log x| + \log\log\frac{1}{\eps}\right)$ arithmetic operations. \qed 
\end{brent} 

\noindent 
While Brent's paper \cite{brent} does not explicitly mention general  
real numbers, he works with a model of floating point number from which 
it is routine to derive the statement above.  

\subsection{Unboundedness and Sign Checks} 
\label{sub:degen} 

\noindent 
Optimizing an $f\!\in\!\cF^{**}_{n,n+1}$ will ultimately 
reduce to checking simple inequalities involving just the coefficients 
of $f$. The optimum will then in fact be either $+\infty$ or the\linebreak 
constant 
term of $f$. Optimizing an $f\!\in\!\cF^{**}_{n,n+2}$ would be as easy were 
it not for two additional difficulties: deciding unboundedness already 
entails checking the sign of a\linebreak 
generalized $\cA$-discriminant, and the optimum 
can be a transcendental function of the coefficients. 

To formalize the harder case, let us now work at the level of 
exponential sums: let us define $\cG_{n,m}$, $\cG^*_{n,m}$, and 
$\cG^{**}_{n,m}$ to be the obvious respective exponential 
$m$-sum analogues of $\cF_{n,m}$, $\cF^*_{n,m}$, and $\cF^{**}_{n,m}$. 
Recall that $\conv \cA$ is the convex hull of $\cA$. 
\begin{thm}
\label{thm:bigger} 
Suppose we write $g\!\in\!\cG^{**}_{n,n+2}$ in the form 
$g(y)\!=\!\sum^{n+2}_{i=1} c_ie^{a_i\cdot y}$ with $\cA\!=\!\{a_1,
\ldots,a_{n+2}\}$. Let us also order the monomials of $f$ so that 
$\cB\!:=\!\{a_1,\ldots,a_{j'}\}$ is the\linebreak 
\scalebox{.92}[1]{unique 
non-degenerate sub-circuit of $\cA$ and let 
$b\!:=\!(b_1,\ldots,b_{n+2})$}\linebreak where
$b_i$ is $(-1)^i$ times the determinant
of the matrix with columns $1\times a_1,\ldots, 
\widehat{1\times a_i}, \ldots,a_{n+2}$ ($\widehat{(\cdot)}$ denoting omission).
Then $\sup_{y\in\Rn} g(y)\!=\!+\infty 
\Longleftrightarrow$ one of the following $2$ conditions hold:
\begin{enumerate}
\item{$c_j\!>\!0$ for some vertex $a_j$ of $\conv\cA$ not equal to $\bO$.} 
\item{$\bO\!\not\in\!\cB$, we can further order the monomials of $f$ so that 
$a_{j'}$ is the unique point of $\cB$ in the relative\linebreak 
\scalebox{.92}[1]{interior of $\cB$, 
$c_{j'}\!>\!0$, and $\prod^{j'}_{i=1} 
\left(\sign(b_{j'})\frac{c_i}{b_i}\right)^{\sign(b_{j'})b_i}\!
<\!1$.}} 
\end{enumerate} 
Finally, if $\sup_{y\in\Rn} g(y)\!=\!\lambda^*\!<\!+\infty$ and 
$a_j\!=\!\bO$, then $\lambda^*\!=\!c_j$, or $\lambda^*$ is 
the unique solution to\\   
\scalebox{.9}[1]{$\left(\sign(b_{j'})\frac{c_{j}-\lambda^*}{b_{j}}
\right)^{\sign(b_{j'})b_{j}} \times \prod 
\limits_{i\in\{1,\ldots,j'\}\setminus\{j\}} 
\left(\sign(b_{j'})\frac{c_i}{b_i} \right)^{\sign(b_{j'})b_i}\!=\!1$} 
\linebreak 
with $(c_j-\lambda^*)b_jb_{j'}\!>0$; where the equation for 
$\lambda^*$ holds iff: 
\begin{enumerate}
\addtocounter{enumi}{2}
\item{$\bO\!\in\!\cB$, we can further order the monomials of $f$ so that 
$a_{j'}$ is the unique point of $\cB$ in the relative interior of $\cB$, 
and $c_{j'}\!>\!0$. }
\end{enumerate}  
\end{thm} 

\noindent
It is easily checked that $c_1b_1b_{j'},\ldots,c_{j'-1}b_{j'-1}b_{j'}\!>\!0$ 
when Conditions 2 or 3 hold. 
While the $3$ cases above may appear complicated, they are easily 
understood from a tropical perspective: our cases above correspond to 
$4$ different families of generalized Viro diagrams that characterize how the 
function $g$ can be bounded from above (or not) on $\Rn$. 
Some representative examples are illustrated below:\\
\begin{picture}(200,220)(0,-120) 
\put(20,0){\epsfig{file=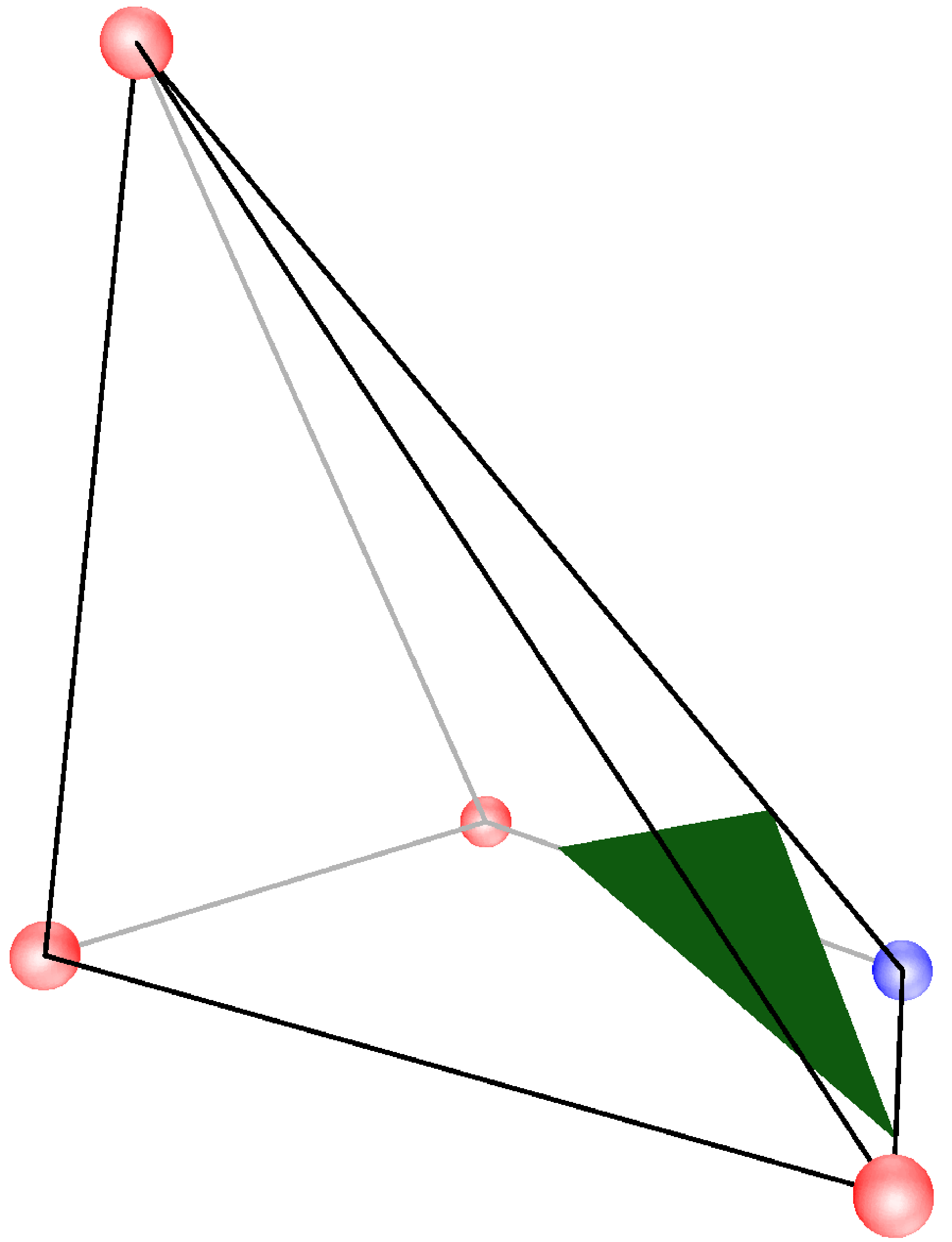,height=1.3in}}
 \put(46,74){Case 1} \put(92,1){{\small $y_1$}} \put(50,36){{\small $y_2$}}
 \put(13,17){{\small $\bO$}} 
\put(130,0){\epsfig{file=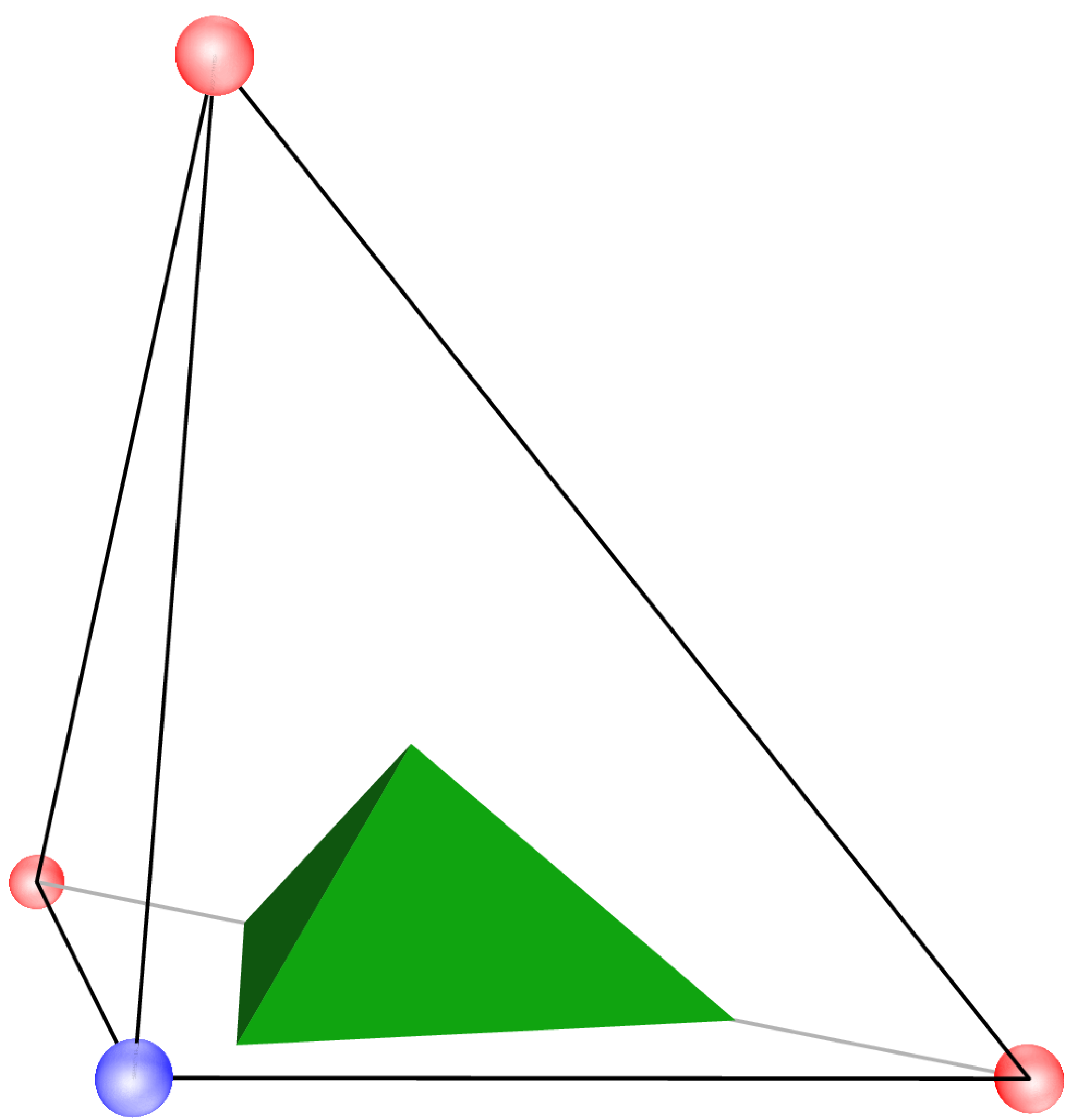,height=1.3in}}
 \put(162,74){Case 2} \put(220,2){{\small $y_1$}} \put(125,25){{\small $y_2$}}
 \put(129,1){{\small $\bO$}}
\put(20,-110){\epsfig{file=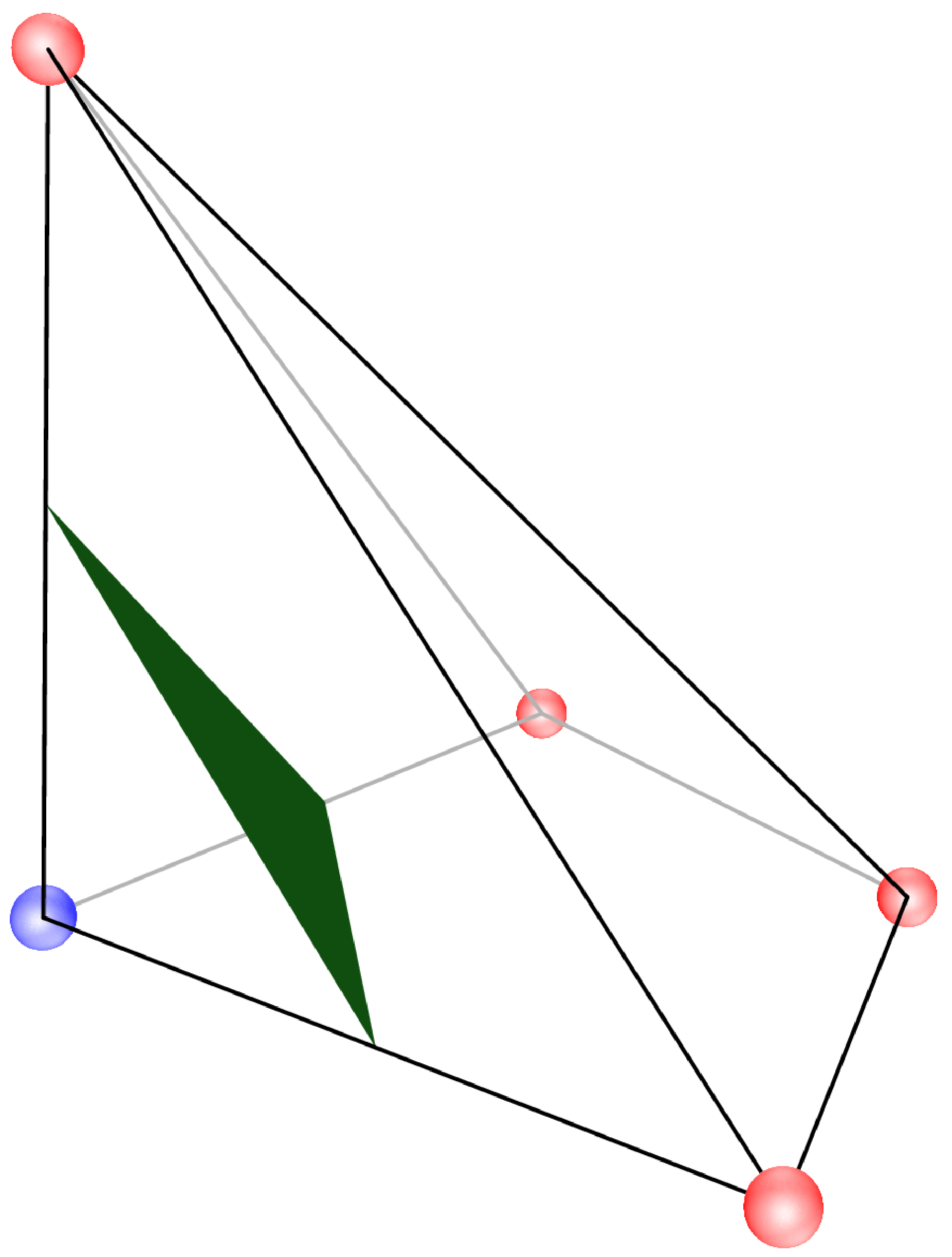,height=1.3in}}
 \put(47,-40){{\bf Not} Case 3} \put(84,-106){{\small $y_1$}} 
 \put(57,-50){($\lambda^*\!<\!+\infty$)} \put(84,-106){{\small $y_1$}} 
 \put(64,-70){{\small $y_2$}} \put(13,-93){{\small $\bO$}}
\put(135,-110){\epsfig{file=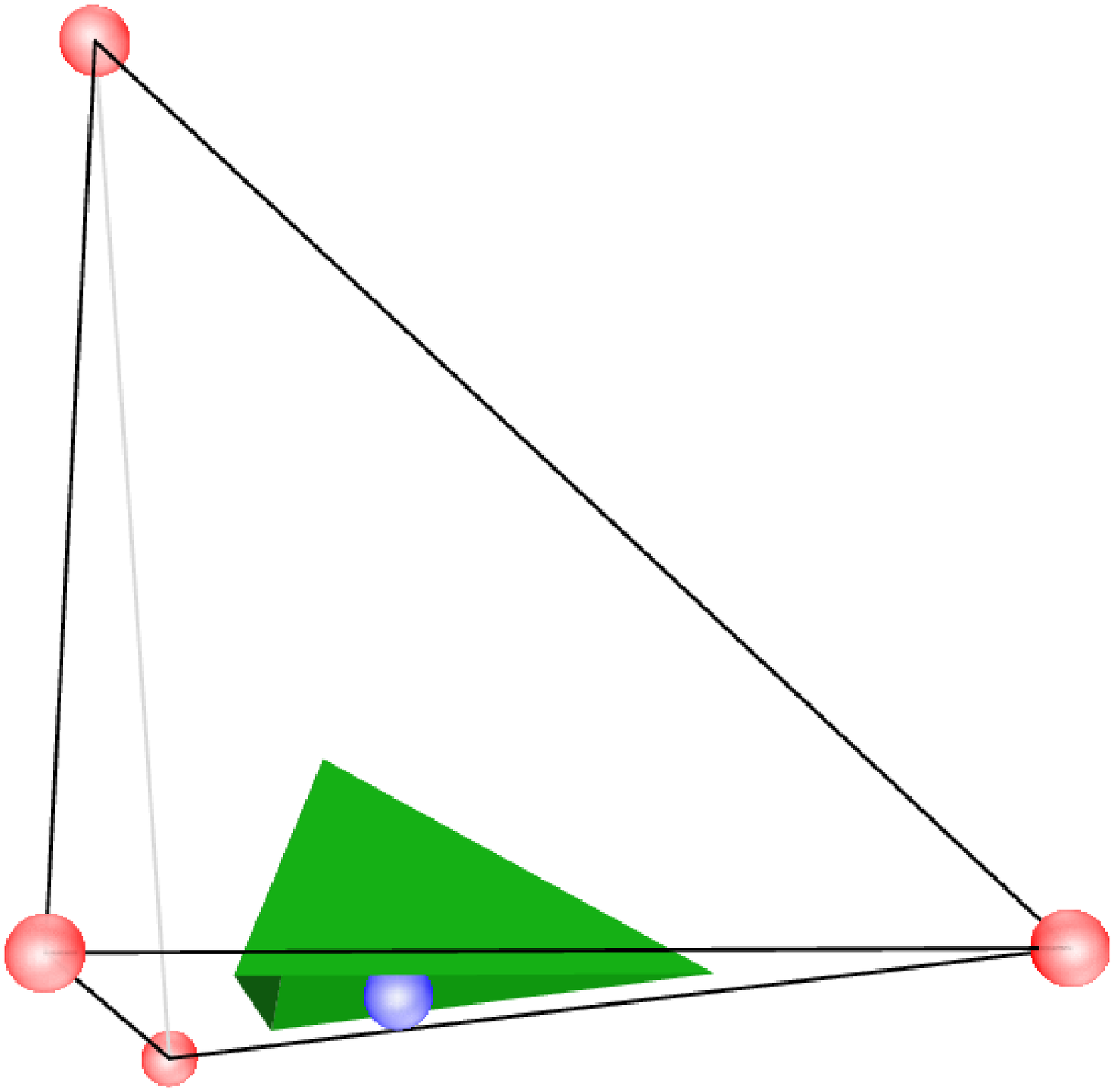,height=1.3in}}
 \put(166,-40){Case 3} \put(229,-102){{\small $y_1$}}
 \put(178,-50){($\lambda^*\!<\!+\infty$)}
 \put(126,-100){{\small $\bO$}} \put(148,-115){{\small $y_2$}}
\end{picture}

\noindent
For example, the first two illustrations are meant to encode the 
fact that there exist directions in the positive quadrant 
along which $g$ increases without bound. Similarly, 
the last $2$ illustrations respectively show cases where 
$g$ either approaches a supremum as some $y_i\longrightarrow-\infty$ 
or has a unique maximum in the real plane. 

\medskip
\noindent 
{\bf Sketch of Proof of Theorem \ref{thm:bigger}:} 
First, we identify the graph of $g$ over $\Rn$ with the real zero set $Z$ of 
$z-g(y)$. Since the supremum of $g$ is unaffected by a linear change of 
variables, we can then assume (analogous to Corollary \ref{cor:can}) that 
$g$ is of the form\\
\mbox{}\hfill $c+e^{y_1}+\cdots+e^{y_\ell}-e^{y_{\ell+1}}-\cdots-e^{y_n}
+c'e^{\alpha\cdot y}$.\hfill\mbox{}\\ 
(Note in particular that a linear change of variables for an 
exponential sum is, modulo applications of $\exp$ and $\log$, the 
same as a monomial change of variables.) Note also that the 
classical Hadamard bound for the determinant guarantees that 
$\log \cC(g)$ increases by at worst a factor of $n$ after our 
change of variables. Let $P$ denote the convex hull of $\{\bO,e_1,
\ldots,e_n,e_{n+1},\alpha\}$.  

Via a minor variation of the {\bf moment map} (see, e.g., \cite{tfulton}) 
one can then give a homeomorphism $\fii : \R^{n+1} \longrightarrow 
\mathrm{Int}(P)$ that extends to a map $\bar{\fii}$ encoding the 
``limits at toric infinity'' of $Z$ 
in terms of data involving $P$. (See also \cite[Sec.\ 6]{tri}.) In particular, 
$\bar{\fii}(Z)$ intersects the facet of $P$ parallel to the 
$y_i$ coordinate hyperplane iff $Z$ contains points with $y_i$ 
coordinates approaching $-\infty$. Similarly, the function $g$ 
is unbounded iff $\bar{\fii}(Z)$ intersects a face of $P$ incident to 
$e_{n+1}$ and some point in $\{e_1,\ldots,e_n,\alpha\}$. This 
correspondence immediately accounts for Condition 1. 

This correspondence also accounts for Condition 2, but in 
a more subtle manner. In particular, $Z$ has topology depending exactly on 
which connected component of the complement of $\nabla_\cA$ contains 
$g$. Thanks to Theorem \ref{THM:DISC}, this can be decided 
by determining the sign of expression involving powers of 
ratios of $c_i$ and $b_i$. In particular, Condition 
2 is nothing more than an appropriate accounting 
of when $\bar{\fii}(Z)$ intersects a face of $P$ incident to
$e_{n+1}$ and some point in $\{e_1,\ldots,e_n,\alpha\}$. 

To conclude, one merely observes that Condition 3 corresponds 
to $\bar{\fii}(Z)$ intersecting a face of $P$ incident to 
$\bO$ and $e_{n+1}$. In particular, the sign conditions merely guarantee that 
$g$ has a unique maximum as some $y_i$ tend to $-\infty$. \qed  

\section{The Proofs of Our Main Results: Theorems \ref{THM:DISC} and 
\ref{THM:BIG}, and Corollaries \ref{COR:RFEAS} and \ref{COR:4}}  
\label{sec:proofs}

\noindent 
We go in increasing order of proof length. 

\subsection{The Proof of Theorem \ref{THM:DISC}}
\label{sub:proof2}

\medskip 
\noindent
{\bf Assertion (1):} It is easily checked
that $Z_\C(f)$ has a degenerate point $\zeta$ iff
\[ \hA \begin{bmatrix}
c_1 e^{a_1\cdot \zeta}\\
\vdots \\
c_{n+2} e^{a_{n+2}\cdot \zeta}\\
\end{bmatrix}
= \begin{bmatrix} 0 \\ \vdots \\ 0 \end{bmatrix}. \]

\noindent
In which case, $(c_1 e^{a_1\cdot y},\ldots,c_{n+2} e^{a_{n+2}\cdot y})^T$ 
must be a generator of the right null space of $\hA$. On the other hand,
by Cramer's Rule, one sees that $(b_1,\ldots,b_{n+2})^T$ is also
a generator of the right null space of $\hA$. In particular,
$\cA$ a non-degenerate circuit implies that $b_i\!\neq\!0$ for all $i$.

We therefore obtain that\\
\mbox{}\hfill $(c_1 e^{a_1\cdot \zeta},\ldots,c_{n+2} 
e^{a_{n+2}\cdot \zeta})\!=\!\alpha
(b_1,\ldots,b_{n+2})$\hfill\mbox{}\\ 
for some $\alpha\!\in\!\Cs$. Dividing
coordinate-wise and taking absolute values, we then obtain\\
$\left(|c_1/b_1|e^{a_1\cdot \re(\zeta)},\ldots,|c_{n+2}/b_{n+2}|
e^{a_{n+2}\cdot \re(\zeta)}\right) \!=\!(|\alpha|,\ldots,|\alpha|)$. Taking 
both sides to the vector power $(b_1,\ldots,b_{n+2})$ we then clearly obtain\\
\scalebox{.84}[1]{$\left(|c_1/b_1|^{b_1}\cdots 
|c_{n+2}/b_{n+2}|^{b_{n+2}}\right)
\left(e^{(b_1 a_1+\cdots b_{n+2}a_{n+2})\cdot \re(\zeta)}\right)\!=\!
|\alpha|^{b_1+\cdots+b_n}$.}\linebreak  
Since $\hA(b_1,\ldots,b_{n+2})^T\!=\!\bO$, we thus obtain
$\prod\limits^{n+2}_{i=1} \left|\frac{c_i}{b_i}\right|^{b_i}\!=\!1$.
Since the last equation is homogeneous in the $c_i$, its zero set
in $\Pro^{n+1}_\C$ actually defines a closed set of $[c_1:\cdots:c_{n+2}]$.
So we obtain the containment for $\nabla_\cA$.

\scalebox{.95}[1]{The assertion on the complexity of computing 
$(b_1,\ldots,b_{n+2})$}\linebreak 
then follows immediately from the classic efficient parallel  
algorithms for linear algebra over $\R$ \cite{csanky}. \qed 

\smallskip
\noindent
{\bf Assertion (2):}   
We can proceed by
almost exactly the same argument as above, using one simple additional
observation:
$e^{a_i\cdot \zeta}\!\in\!\R_+$ for all $i$ when $\zeta\!\in\!\R$. So then,
we can replace our use of absolute value by a sign factor, so that
all real powers are well-defined. In particular, we immediately obtain
the ``$\Longleftarrow$'' direction of our desired equivalence.

To obtain the ``$\Longrightarrow$'' direction, note that when\\
\mbox{}\hfill 
$Z_\R\!\left(\sum^{n+2}_{i=1}
c_i e^{a_i\cdot y}\right)$ \hfill \mbox{}\\ 
has a degeneracy $\zeta$, we directly
obtain $e^{a_i\cdot \zeta}\!=\!
\sign(b_1c_1)b_i/c_i$ for all $i$ (and the constancy of $\sign(b_ic_i)$
in particular). We thus obtain the system of equations\\
\mbox{}\hfill
$\left(e^{(a_2-a_1)\cdot \zeta},\ldots,e^{(a_{n+1}-a_1)\cdot \zeta}\right)
=\left(\frac{b_2c_1}{b_1c_2},\ldots,\frac{b_{n+1}c_1}{b_1c_{n+1}}\right)$,
\hfill\mbox{}\\
and $a_2-a_1,\ldots,a_{n+1}-a_1$ are linearly independent since $\cA$ is a
circuit. So, employing Proposition \ref{prop:monochange}, we can easily solve
the preceding system for $\zeta$ by taking the logs of the coordinates of
$\left(\frac{b_2c_1}{b_1c_2},\ldots,
\frac{b_{n+1}c_1}{b_1c_{n+1}}\right)^{[a_2-a_1,\ldots,a_{n+1}-a_1]^{-1}}$.
\qed 

\subsection{Proving Corollary \ref{COR:RFEAS} and Theorem \ref{THM:BIG} } 
\label{sub:thresh} 

\medskip 
\noindent 
{\bf Corollary \ref{COR:RFEAS} and Assertion (2) of Theorem 
\ref{THM:BIG}:}\linebreak   
\scalebox{.96}[1]{Since our underlying family of putative hard problems 
shrinks}\linebreak as $\delta$ decrease, it clearly suffices to prove the 
case $\delta\!<\!1$. 
So let assume henceforth that $\delta\!<\!1$. Let us also define 
$\qsat_\R$ to be the problem of deciding whether an input  
{\bf quartic} polynomial $f\!\in\!\bigcup_{n\in\N} \R[x_1,\ldots,x_n]$ 
has a real root or not. $\qsat_\R$ (referred to as $4$-FEAS in \cite{bcss}) 
is one of the fundamental $\np_\R$-complete problems (see Chapter 
4 of \cite{bcss}). 

That $\opt\!\in\!\np_\R$ follows immediately from 
the definition of $\np_\R$. So it suffices to prove that\\ 
\mbox{}\hfill
$\opt\!\left(
\bigcup\limits_{\substack{n\in\N\\ 0<\delta'<\delta}}
\cF^{**}_{n,n+n^{\delta'}}\cap\R[x_1,\ldots,x_n]
\right)$\hfill\mbox{}\\ 
is $\np_\R$-hard. We will do this by giving an explicit 
reduction of $\qsat_\R$ to\\   
\mbox{}\hfill $\opt\!\left(
\bigcup\limits_{\substack{n\in\N\\ 0<\delta'<\delta}}
\cF^{**}_{n,n+n^{\delta'}}\cap\R[x_1,\ldots,x_n]\right)$,\hfill\mbox{}\\ 
passing through $\feas_+
\text{\raisebox{-.3cm}{\scalebox{1}[2.3]{$($}}} 
\bigcup\limits_{\substack{n\in\N\\ 0<\delta'<\delta}}
\cF^{**}_{n,n+n^{\delta'}}\cap\R[x_1,\ldots,x_n]
\text{\raisebox{-.3cm}{\scalebox{1}[2.3]{$)$}}}$ 

\vspace{-.4cm} 
\noindent 
along the way. 

\medskip
To do so, let $f$ denote any $\qsat_\R$ instance, involving, say, 
$n$ variables. Clearly, $f$ has no more than \scalebox{.7}[.7]
{$\begin{pmatrix} n+4\\ 4\end{pmatrix}$}\linebreak 
monomial terms. 
Letting $\qsat_+$ denote the natural\linebreak variant of $\qsat_\R$ where one 
instead asks if $f$ has a root in $\Rn_+$, we will first need to 
show that $\qsat_+$ is $\np_\R$-hard as an intermediate step. 
This is easy, via the introduction of slack variables: using $2n$ new 
variables $\left\{x^\pm_i\right\}^n_{i=1}$ 
and\linebreak 
forming the polynomial $f^\pm(x^\pm)\!:=\!f\!\left(x^+_1-x^-_1,\ldots,
x^+_n-x^-_n\right)$, it is clear that $f$ has a root in $\Rn$ iff 
$f^\pm$ has a root in $\R^{2n}_+$. Furthermore, we easily see that\\ 
\mbox{}\hfill $\size(f^\pm)\!=\!(16+o(1))\size(f)$.\hfill\mbox{}\\ 
So $\qsat_+$ is 
$\np_\R$-hard. We also observe that we may\linebreak 
restrict the inputs 
to quartic polynomials with full-\linebreak 
dimensional Newton polytope, since 
the original proof for the $\np_\R$-hardness of $\qsat_\R$ actually 
involves polynomials having nonzero constant terms and nonzero 
$x^4_i$ terms for all $i$ \cite{bcss}. 

So now let $f$ be any $\qsat_+$ instance with, say, $n$\linebreak variables. 
Let us also define, for any $M\!\in\!\N$, the polynomial 
$t_M(z)\!:=\!1+z^{M+1}_1+\cdots+z^{M+1}_M-(M+1)z_1\cdots z_M$. 
One can then check via the Arithmetic-Geometric Inequality 
\cite{hlp} that $t_M$ is nonnegative on $\R^M_+$, with a unique 
root at $z\!=\!(1,\ldots,1)$. 
Note also that $f^2$ has no more than \scalebox{.7}[.7]{$\begin{pmatrix}
n+4\\ 4\end{pmatrix}^2$} monomial terms. 
Forming the polynomial $F(x,z)\!:=\!f(x)^2+t_M(z)$ 
with $M\!:=\!\left\lceil \begin{pmatrix} n+4\\ 4\end{pmatrix}^{2/\delta}
\right\rceil$, we see that $f$ has a root in $\Rn_+$ iff $F$ has a root in 
$\R^{n+M}_+$. It is also easily checked that $F\!\in\!\cF^{**}_{N,N+k}$ with 
$k\!\leq\!N^{\delta'}$, where $N\!:=\!n+M$ and $0\!<\!\delta'\!\leq\!\delta$. 
In particular,\\ 
\mbox{}\hfill 
$k\!<\!\begin{pmatrix} n+4\\ 4\end{pmatrix}^2\!\leq\!\left\lceil\begin{pmatrix}
n+4\\ 4\end{pmatrix}^{2/\delta}\right\rceil^\delta\!=\!M^\delta\!<\!
(n+M)^\delta$.\hfill\mbox{}\\ 
So we must now have that\\
\mbox{}\hfill  
$\feas_+
\text{\raisebox{-.3cm}{\scalebox{1}[2.3]{$($}}} 
\bigcup_{\substack{n\in \N\\ 0<\delta'<\delta}} 
\cF^{**}_{n,n+n^\delta}\cap\R[x_1,\ldots,x_n]
\text{\raisebox{-.3cm}{\scalebox{1}[2.3]{$)$}}} 
$\hfill\mbox{}\\ 
is $\np_\R$-hard. (A small digression allows us to succinctly prove that\\ 
\mbox{}\hfill $\feas_\R 
\text{\raisebox{-.3cm}{\scalebox{1}[2.3]{$($}}} 
\bigcup_{\substack{n\in \N\\ 0<\delta'<\delta}}\cF^{**}_{n,n+n^\delta}
\cap\R[x_1,\ldots,x_n]
\text{\raisebox{-.3cm}{\scalebox{1}[2.3]{$)$}}} 
$\hfill\mbox{}\\ 
is $\np_\R$-hard as well: we simply repeat the argument from the last 
paragraph, but use $\qsat_\R$ in place of $\qsat_+$, and define 
$F(x,z)\!:=\!f(x)^2+t_M(z^2_1,\ldots,z^2_M)$ instead.)  

To conclude, note that $F(x,z)$ is nonnegative on $\Rn_+$. 
So by checking whether $-F$ has supremum $\geq\!0$ in $\Rn_+$, we can 
decide if $F$ has a root in $\Rn_+$. In other words,\\ 
\mbox{}\hfill 
$\opt\!\left(\bigcup\limits_{\substack{n\in\N\\ 0<\delta'<\delta}}
\cF^{**}_{n,n+n^{\delta'}}\cap\R[x_1,\ldots,x_n]\right)$\hfill \mbox{}\\ 
must be $\np_\R$-hard has well. So we are done. \qed 

\medskip 
\noindent 
{\bf Assertion (0) of Theorem \ref{THM:BIG}:} Letting $(f,\eps)$ denote 
any instance of $\oopt\!\left(\bigcup_{n\in\N} \cF^{**}_{n,n+1}\right)$, 
first note that via\linebreak Corollary \ref{cor:can} 
we can assume that\\ 
\mbox{}\hfill $f(x)\!=\!c_1+x_1+\cdots+x_\ell-x_{\ell+1}-\cdots 
-x_n$,\hfill\mbox{}\\ 
after a computation in $\nc^1_\R$. Clearly then, $f$ has an\linebreak 
unbounded 
supremum iff $\ell\!\geq\!1$. Also, if $\ell\!=\!0$, then 
the\linebreak supremum of $f$ is exactly $c_1$. So 
$\oopt\!\left(\bigcup_{n\in\N} 
\cF^{**}_{n,n+1}\right)\!\in\!\nc^1_\R$. That 
$\opt\!\left(\bigcup_{n\in\N} \cF^{**}_{n,n+1}\right)\!\in\!\nc^1_\R$ is 
obvious as well: after checking the signs of the $c_i$, we make merely 
decide the sign of $c_1-\lambda$. \qed 

\begin{rem}
\label{rem:nc1}
Note that checking whether a given $f\!\in\!\cF_{n,n+1}$ 
lies in $\cF^*_{n,n+1}$ can be done within $\nc^2$: 
one\linebreak simply finds $d\!=\!\dim \supp(f)$ in $\nc^2$ by 
computing the rank of the matrix whose columns are $a_2-a_1,
\ldots,a_m-a_1$ (via the parallel algorithm of Csanky \cite{csanky}), 
and then checks whether $d\!=\!n$. \dia 
\end{rem}

\noindent 
{\bf Assertion (1):} We will first derive the $\hhpt$ result. 
Let us assume $f\!\in\!\cF^{**}_{n,n+2}$ and observe the following 
algorithm: 
\newpage 

\begin{algor} 
\label{algor:fsup} 
\mbox{}\\
{\bf Input:} A coefficient vector $c\!:=\!(c_1,
\ldots,c_{n+2})$, a (possibly degenerate) circuit
$\cA\!=\!\{a_1,\ldots,a_{n+2}\}$ of cardinality $n+2$, and 
a precision parameter $\eps\!>\!0$. \\
{\bf Output:} A pair 
\mbox{}\hfill $(\bar{x},\bar{\lambda})\!\in\!(\R_+\cup\{0,+\infty\})^n\times
(\R\cup\{+\infty\})$\hfill \mbox{}\\
with $\bar{x}\!=\!(\bar{x}_1,\ldots,\bar{x}_n)$ (resp.\
$\bar{\lambda}$) an $\hhpt$ for $x^*$ (resp.\ $\lambda^{*}$)
where $f(x)\!:=\!\sum^{n+2}_{i=1}c_ix^{a_i}$ and 
$\lambda^*\!:=\!\sup_{x\in\Rn_+} f\!=\!\lim_{x\rightarrow x^*}f(x)$
for some $x^{*}\!=\!(x^*_1,\ldots,x^*_n)\!\in\!(\R_+\cup\{0,+\infty\})^n$. 

\medskip 
\noindent
{\bf Description:}

\vspace{-.2cm}
\begin{enumerate}
\item{If $c_i\!>\!0$ for some $i$ with $a_i\!\neq\!\bO$ a vertex of 
$\conv \cA$ then output\\ ``{\tt $f$ tends to $+\infty$ along a 
curve of the form\\ $\{ct^{a_i}\}_{t\rightarrow +\infty}$}''\\ 
and {\tt STOP}.} 
\item{Let $b\!:=\!(b_1,\ldots,b_{n+2})$ where
$b_j$ is $(-1)^j$ times the\linebreak determinant
of the matrix with columns $1\times a_1,\ldots,$\linebreak 
$\widehat{1\times a_j}, \ldots,a_{n+2}$ ($\widehat{(\cdot)}$ denoting omission).
If $b$ or $-b$ has a unique negative
coordinate $b_{j'}$, and $c_{j'}$ is the unique negative coordinate of $c$, 
then do the following:
\begin{enumerate}
\item{Replace $b$ by $-\sign(b_{j'})b$ 
and then reorder $b$, $c$, and $\cA$ by the same 
permutation so that $b_{j'}\!<\!0$ and [$b_i\!>\!0$ iff $i\!<\!j'$]. }
\item{If $a_i\!\neq\!\bO$ for all $i\!\in\!\{1,\ldots,j'\}$ and\\
\mbox{}\hfill  
$\prod^{j'}_{i=1}
\left(\sign(b_{j'})\frac{c_i}{b_i}\right)^{\sign(b_{j'})b_i}\!<\!1$ 
\hfill\mbox{}\\
then output\\ ``{\tt $f$ tends to $+\infty$ along a
curve of the form\\ $\{ct^{a_{j'}}\}_{t\rightarrow +\infty}$}''\\
and {\tt STOP}.} 
\item{If $a_j\!=\!\bO$ for some $j\!\in\!\{1,\ldots,j'\}$ then 
output\\ 
``{\tt $f(z)$ tends to a supremum of $\bar{\lambda}$ as 
$z$ tends\\ to the point $\bar{x}$ on the $(j'-2)$-dimensional\\ 
sub-orbit corresponding to $\{a_1,\ldots,a_{j'}\}$.}'',\\ 
where $x\!\in\!\R^{j'-2}_+$ is the unique solution to the\linebreak 
binomial system\\ 
\mbox{}\hfill
$\left(x^{a_2-a_1},\ldots,x^{a_{j'-1}-a_1}\right)
=\left(\frac{b_2c_1}{b_1c_2},\ldots,\frac{b_{n+1}c_1}{b_1c_{n+1}}\right)$,
\hfill\mbox{}\\
$\bar{x}$ is a $(1+\eps)$-factor approximation\footnote{We compute $\bar{x}$ 
and $\bar{\lambda}$ via Proposition \ref{prop:monochange} and the 
Brent-Salamin Theorem.} of $x$, $\lambda$ is the unique solution of\\
\mbox{}\hspace{-1.3cm}
\scalebox{.8}[1]{$\left(\sign(b_{j'})\frac{c_{j}-\lambda}{b_{j}}
\right)^{\sign(b_{j'})b_{j}} \times \prod
\limits_{i\in\{1,\ldots,j'\}\setminus\{j\}}
\left(\sign(b_{j'})\frac{c_i}{b_i} \right)^{\sign(b_{j'})b_i}\!=\!1$}
\linebreak
with $(c_j-\lambda)b_jb_{j'}\!>0$, and 
$\bar{\lambda}$ is 
a $(1+\eps)$-factor approximation$^8$ of $\lambda$; 
and {\tt STOP}.}
\end{enumerate} }
\item{Output\\ 
``{\tt $f$ approaches a supremum of $c_j$ as 
all $x^{a_i}$ with\\ $a_i$ incident to $a_j$ approach $0$.}'',\\ 
where $a_j\!=\!\bO$, and {\tt STOP}.}
\end{enumerate}
\end{algor} 

Our proof then reduces to proving correctness and a suitable complexity 
bound for Algorithm \ref{algor:fsup}. In particular, 
correctness follows immediately from Theorem \ref{thm:bigger}. So we now 
focus on a complexity analysis.  

Steps 1 and 3 can clearly be done within $\nc^1_\R$, so let us 
focus on Step 2. 

For Step 2, the dominant complexity comes from Part (b). 
(Part (a) can clearly be done in $\nc^1_\R$, and Part (c) can clearly 
be done in $\nc^2_\R$ via Csanky's method \cite{csanky}.) 
The latter can be done by taking the logarithm 
of each term, thus reducing to checking the sign of a linear combination of 
logarithms of positive real numbers. 
So the arithmetic complexity of our algorithm is  
$O\!\left(\log \cC(f)+\log\log\frac{1}{\eps}\right)$ and we thus 
obtain our $\hhpt$ result. 

The proof that 
$\opt\!\left(\bigcup_{n\in\N}\cF^{**}_{n,n+2}\right)\!\in\!\pp_\R$ 
is almost completely identical. \qed 

\medskip
Note that just as in Remark \ref{rem:nc1}, checking whether a given 
$f\!\in\!\cF_{n,n+2}$ lies in $\cF^*_{n,n+2}$ can be done within $\nc^2$ 
by \linebreak 
computing $d\!=\!\dim \supp(f)$ efficiently. Moreover, 
deciding whether a circuit is degenerate (and 
extracting $\cB$ from $\cA$ when $\cA$ is degenerate) can be done in $\nc^2$ 
as well since this is ultimately the evaluation of $n+2$ determinants. 

\subsection{The Proof of Corollary \ref{COR:4}} 
\label{sub:4} 

\noindent 
{\bf Assertion (0):} Since the roots of $f$ in $\R_+$ are 
unchanged under multiplication by monomials, we can clearly 
assume $f\!\in\!\cF^{**}_{1,3}\cap\R[x_1]$. Moreover, via 
the classical Cauchy bounds on the size of roots of 
polynomials, it is easy to show that the log of any root of $f$ 
is $O(\log\cC(f))$. We can then invoke Theorem 1 of \cite{rojasye} 
to obtain our desired $\hhpt$ as follows: 
If $D\!:=\!\deg(f)$, 
\cite[Theorem 1]{rojasye} tells us that we can count exactly the number of 
positive roots of 
$f$ using $O(\log^2 D)$ arithmetic operations, and $\eps$-approximate 
all the roots of $f$ in $(0,R)$ within 
$O\!\left((\log D)\log\left(D\log\frac{R}{\eps}\right)\right)$ arithmetic 
operations. Since we can take $\log R\!=\!O(\log \cC(f))$ via our 
root bound observed above, we are done. \qed 

\medskip 
\noindent 
{\bf Assertion (1):} 
Writing any $f\!\in\!\cF^{**}_{1,4}\cap\R[x_1]$ as 
$f(x)\!=\!c_1+c_2x^{a_2}+c_3x^{a_3}+c_4x^{a_4}$ with 
$a_2\!<\!a_3\!<\!a_4$, note 
that $f$ has unbounded supremum on $\R_+$ iff $c_4\!>\!0$
So let us assume $c_4\!<\!0$.  

Clearly then, the supremum of $f$ is attained either at a 
critical point in $\R_+$ or at $0$.
But then, any positive critical point is a positive root of a trinomial, 
and by Assertion (0), such critical points must admit an 
$\hhpt$. Similarly, since $f$ is a tetranomial (and thus 
evaluable within $O(\log \deg(f))$ arithmetic operations), 
we can efficiently approximate (as well as 
efficiently check inequalities involving) $\sup_{x\in\R_+} f$. 
So we are done. \qed  

\section*{Acknowledgements} 
We thank Peter B\"urgisser, Felipe Cucker, Johan Hastad, and 
Gregorio Malajovich for discussions on complexity classes over $\R$. 
The second author would also like to thank MSRI and the Wenner Gren 
Foundation for their support during the completion of this paper. 
In particular, special thanks go to Mikael Passare and 
Boris Shapiro of Stockholm University for their hospitality.

\bibliographystyle{abbrv}

\end{document}